\documentclass[12pt,leqno]{article}
\usepackage{amsthm,amsbsy,amsfonts,amssymb,amsmath,amscd}
\usepackage{latexsym}
\usepackage{fontenc}
\usepackage[mathscr]{euscript}
\usepackage{graphics}

\DeclareMathAlphabet{\scrb}{U}{eus}{b}{n}
\DeclareMathAlphabet{\cur}{U}{eur}{b}{n}

\newcommand{\bfE}{{\mathbf E}}
\newcommand{\bfe}{{\mathbf e}}

\newcommand{\ScrE}{{\mathscr{E}}} 
\newcommand{\ScrH}{{\mathscr{H}}} 
\newcommand{\ScrL}{{\mathscr{L}}} 
\newcommand{\ScrN}{{\mathscr{N}}} 
\newcommand{\ScrO}{{\mathscr{O}}} 
\newcommand{\ScrR}{{\mathscr{R}}} 
\newcommand{\ScrT}{{\mathscr{T}}}
\newcommand{\ScrW}{{\mathscr{W}}}

\newcommand{\bfscrE}{{\scrb{E}}}

\font\teneurm=eurm10 \font\seveneurm=eurm7 
\font\fiveeurm=eurm5 
\newfam\eurmfam 
\textfont\eurmfam=\teneurm 
\scriptfont\eurmfam=\seveneurm 
\scriptscriptfont\eurmfam=\fiveeurm

\theoremstyle{plain}
\newtheorem{thm}{Theorem}

\newtheorem{lem}{Lemma}

\theoremstyle{remark}
\newtheorem*{Rem}{\bf Remark}

\newcommand{\Le}{{{\mathchoice{\,{\scriptstyle\le}\,}
  {\,{\scriptstyle\le}\,}
  {\,{\scriptscriptstyle\le}\,}{\,{\scriptscriptstyle\le}\,}}}}
\newcommand{\Ge}{{{\mathchoice{\,{\scriptstyle\ge}\,}
  {\,{\scriptstyle\ge}\,}
  {\,{\scriptscriptstyle\ge}\,}{\,{\scriptscriptstyle\ge}\,}}}}
\newcommand{\one}{\hbox{1}\kern-4.25pt \vrule height 6.50pt
depth 0pt width.20pt\,\,}

\font\titlefont=cmss10 scaled\magstep3
\font\sectionfont=cmss10 scaled\magstep3
\font\subsectionfont=cmss10 scaled\magstep2
\font\itemfont=cmss10 scaled\magstep1 

\newcommand{\plus}{{\scriptstyle +}}

\newcommand{\util}{\tilde{u}}
\newcommand{\Ftil}{\widetilde{F}}

\newcommand{\vhat}{\hat{v}}
\newcommand{\xhat}{\hat{x}}
\newcommand{\zhat}{\hat{z}}

\newcommand{\xdot}{\Dot{x}}
\newcommand{\ydot}{\Dot{y}}
\newcommand{\zdot}{\Dot{z}}

\newcommand{\re}{\text{\rm Re}\,}

\title{\titlefont Existence of Chaos in Evolution Equations}
\author{Yanguang (Charles) Li \thanks{This work is supported by 
a Guggenheim Fellowship, and an AMS Centennial Fellowship. Current 
Address: Department of Mathematics, University of Missouri, Columbia, 
MO 65211. Email: cli@math.missouri.edu} \\ School of Mathematics \\ 
Institute for Advanced Study \\ Princeton, NJ 08540}

\begin{document}

\setlength{\unitlength}{1in}

\maketitle

\begin{abstract}
For a general evolution equation with a Silnikov homoclinic
orbit, Smale horseshoes are constructed with the tools of
\cite{Li99} and in the same way as in \cite{Li99}.
The linear part of the evolution equation has a finite
number of unstable modes.
For evolution equations with infinitely many linearly
unstable modes, the problem is still open.

MSC: 35, 37.
\end{abstract}

\section{\sectionfont Introduction}
In recent years, the existence of chaos in partial
differential equations has been established 
\cite{LM94}, \cite{LMSW96}, \cite{Li99}.
These studies set up a scheme for attacking problems 
on Chaos in PDEs.
The type of chaos studied in these works is the so-called
homoclinic chaos generated in a neighborhood of a homoclinic
orbit.
Two types of homoclinic orbits have been studied.
One type is the so-called transversal homoclinic orbit
\cite{Li02}.
The other type is the so-called Silnikov homoclinic orbit,
which is non-transversal \cite{LM94}, \cite{LMSW96},
\cite{Li99}.
For lower dimensional and general finite dimensional
systems, Silnikov studied the symbolic dynamics
structures in the neighborhoods of such homoclinic orbits
\cite{Sil65}, \cite{Sil67a}, \cite{Sil67b}, \cite{Sil70}.
In \cite{LW97} and \cite{Li99}, we have developed a
different construction of Smale horseshoes in the
neighborhood of a Silnikov homoclinic orbit.
The advantages of our constructions have been fully
addressed in \cite{LW97}, \cite{Li99}.
In this note, we generalize the construction in \cite{Li99}
to more general evolution equations with finitely many
linearly unstable modes.
For evolution equations with infinitely many linearly
unstable modes, the problem is still open.

\section{\sectionfont The Set-Up}
Consider the evolution equation
\begin{equation}
\partial_tu=\ScrL u+\ScrN(u)\,\,,
\label{ee}
\end{equation}
where $\ScrL$ is a linear operator which is constant in
time, and $\ScrN$ is the nonlinear term.

Following are the assumptions for the setup:
\begin{description}
\item[\itemfont{(A1)}]
$u=0$ is a saddle, the linear operator $\ScrL$ has only
point spectrum as follows,
$$
\sigma(\ScrL)=\left\{-\alpha\pm i\beta,\ \gamma,
\ \lambda_j^\pm,\,\,
j\in S^\pm\subset Z^{\plus}\right\}\,\,,
$$
where $\alpha>0$, $\beta>0$, $\gamma>0$, $\re\{\lambda_j^\pm \}
\gtrless 0$; the number of elements in $S^{\plus}$ is
finite, denoted by $N$, $\alpha<\gamma$, and
\[
\alpha<\lambda^{-} =\inf\limits_{j\in
S^-}\{-\re\{\lambda_j^-\}\}\ ,\quad
\gamma<\lambda^{\plus} =\inf\limits_{j\in S^{\plus}}
\{\re\{\lambda_j^{\plus}\}\}\,\,.
\]

\item[\itemfont{(A2)}]
The evolution equation (\ref{ee}) is globally well posed in a Hilbert
space $\ScrH$, that is, there exists a unique solution to
the Cauchy problem of (\ref{ee}), $u(t,u_0)\in
C^0[(-\infty,\infty),\ScrH]$, $u(0,u_0)=u_0$.
Moreover, we assume the regularity condition on initial data
that for any $t\in(-\infty,\infty)$, the evolution operator
$F(u_0)^t=u(t,u_0)$ is $C^n$ in $u_0$ for some
$n\Ge 2$.

\item[\itemfont{(A3)}]
With respect to the saddle $u=0$, the evolution operator
$F^t(u)$ admits a $C^2$ smooth linearization, i.e.
there exists a $C^2$ diffeomorphism
$\ScrR\colon\,\ScrH\to\ScrH$, such that in terms of the new
variable $\util=\ScrR\,u$, the evolution equation (\ref{ee})
is transformed into the linear form
$$
\partial_t\util=\ScrL\util
$$
in a neighborhood of $u=0$.
The conjugated evolution operator
$$
\Ftil^t=\ScrR\,F^t\ScrR^{-1}
$$
is still $C^0$ in time $t$ and $C^2$ in $\util$.
In fact, we assume that the evolution equation (\ref{ee})
takes the normal form:
$$
\begin{cases}
\xdot=-\alpha x-\beta  y+G_x(x,y,z,v^{\pm}), &\\
\ydot=\beta x-\alpha y+G_y(x,y,z,v^{\pm}), &\\
\zdot=\gamma z+G_z(x,y,z,v^{\pm}), &\\
\partial_tv^{\pm}=
  L^{\pm}v^{\pm}+G_{v^{\pm}}(x,y,z,v^{\pm}), &
\end{cases}
$$
where $G=0$ in a neighborhood $\Omega$ of $0$, 
$v^{\plus}=(v_1^{\plus},\dotsc,v_N^{\plus})'$, $N$ is the
number of elements in $S^{\plus}$; $x$, $y$, $z$, and
$v_j^{\plus}{'}$s are real variables, and
\[
\Vert e^{L^{\plus}t}\Vert \Le c^{\plus}e^{\lambda^{\plus} t}\,\,,
\ \ \text{as $t\to-\infty$}\,\,,\quad 
\Vert e^{L^- t}\Vert \Le c^- e^{-\lambda^- t}\,\,,
\ \ \text{as $t\to+\infty$}\,\,.
\]
For references on such linearization results, see for
example \cite{Nik86} etc.

\item[\itemfont{(A4)}]
There exists a Silnikov homoclinic orbit $h(t)$ asymptotic
to $0$.
As $t\to+\infty$, $h$ is tangent to the $(x,y)$-plane at
$0$, and as $t\to-\infty$, $h$ is tangent to the $z$-axis
(without loss of generality, positive $z$-axis) at $0$.
The stable and unstable manifolds of $0$ are $C^2$ smooth,
and
$$
\dim\left\{\ScrT_v W^u\cap \ScrT_v W^s\right\}=1\,\,,
$$
where $v\in h(t)$, $\ScrT_vW^u$ is the tangent space of the
unstable manifold of $0$ at the point $v$ on the homoclinic
orbit $h(t)$, similarly for $\ScrT_v W^s$.
\end{description}

\begin{Rem}
Proving the existence of Silnikov homoclinic orbits in partial 
differential equations is a rather nontrivial question. So far,
this has been done for perturbed nonlinear Schr\"odinger equation
\cite{LMSW96} \cite{Li02a}, perturbed vector nonlinear Schr\"odinger 
equation \cite{Li02d}, and perturbed discrete nonlinear 
Schr\"odinger equation \cite{LM97}. The above and later assumptions 
have been either verified or discussed for these equations in 
\cite{Li99} \cite{Li02b} \cite{LW97}. The perturbed 
Davey-Stewartson II equation has been studied along this direction
\cite{Li00a} \cite{Li02c}. Unfortunately, existence of Silnikov 
homoclinic orbits has not been proved due to some technical difficulty
\cite{Li02c}. I would like to comment on equations that have the potential 
of being casted into the above setup: 1. perturbations of the modified 
KdV equation
\[
u_t +6 u^2 u_x +u_{xxx}=0\ ,
\]
2. perturbations of the derivative nonlinear Schr\"odinger equation
\cite{SD99} \cite{GDY01} \cite{MOMT76} \cite{WKI79}
\[
iu_t=u_{xx}+i\alpha (|u|^2u)_x +2|u|^2u\ , \quad \alpha > 0 \ ,
\]
and 3. perturbations of the derivative nonlinear Schr\"odinger equation
\[
iu_t=u_{xx}-i\alpha u^2\bar{u}_x +2|u|^2u +2\alpha^2|u|^4u \ .
\]
All of the above three equations are integrable systems.
\end{Rem}

\section{\sectionfont The Construction of Smale Horseshoes}

\subsection{\subsectionfont Definitions}

\begin{description}
\item[\itemfont{Definition.}]
The Poincar\'e section $\sum_0$ is defined by the
constraints:
\begin{eqnarray*}
& & y=0,\,\,\eta\exp\{-2\pi\alpha/\beta\}<x<\eta\,\,,\\
& & 0<z<\eta,\,\,\Vert v^{\pm}\Vert<\eta\,\,;
\end{eqnarray*}
where $\eta$ is a sufficiently small constant so that
$\sum_0$ is included in the neighborhood $\Omega$ of $0$
where the dynamics is given by the linear system.

\item[\itemfont{Definition.}]
The auxiliary section $\sum_0^{\plus}$ is defined by the
constraints:
\begin{eqnarray*}
& & y=0,\,\,\eta\,\exp\{-2\pi \alpha/\beta\}<x<\eta\,\,,\\
& & -\eta<z<\eta,\,\,\Vert v^{\pm}\Vert<\eta\,\,.
\end{eqnarray*}
The homoclinic orbit $h$  intersects the $(z=0)$-boundary of
$\sum_0$ at $w^{(\plus)}$ with coordinates denoted by
$$
x=x_*,\,\,
y=0,\,\,
z=0,\,\,
v^{\plus}=0,\,\,
v^-=v_*^-\,\,.
$$
There exists $T>0$ such that the point 
$w^{(-)}=F^{-T}(w^{(\plus)})$ on $h$ (where $F^t$ is the
evolution operator) has the $z$ coordinate equal to $\eta$.
Denote the coordinates of $w^{(-)}$ by
$$
x=y=v^-=0,\,\,
z=\eta,\,\,
v^{\plus}=v_*^{\plus}\,\,.
$$

\item[\itemfont{Definition.}]
The Poincar\'e section $\sum_1$ is defined as:
$$
\sum\nolimits_1=\left(F^{-T}\circ\sum\nolimits_0^{\plus}\right)\cap
\Omega\,\,.
$$

\item[\itemfont{Definition.}]
The map $P_0^1$ from $\sum_0$ to $\sum_1$ is defined as:
\begin{equation*}
\begin{split}
P_0^1\colon\,\, 
&U_0\subset\sum\nolimits_0\longmapsto\sum\nolimits_1\,\,,\\
&\forall\,w\in U_0,\,\,P_0^1(w)=F^{t_*}(w)\in \sum\nolimits_1\,\,,
\end{split}
\end{equation*}
where $t_*=t_*(w)>0$ is the smallest time $t$ such that
$F^t(w)\in\sum_1$.

The map from $\sum_1$ to $\overline{\sum}_0$
($=\sum_0\ \bigcup\ \partial\sum_0$) is defined as:
\begin{equation*}
\begin{split}
P_1^0\colon\,\, &U_1\subset\sum\nolimits_1\longmapsto
  \overline{\sum}_0\,\,,\\
&\forall\, w\in\sum\nolimits_1,\,\,P_1^0(w)=F^T(w)\in
\overline{\sum}_0\,\,.
\end{split}
\end{equation*}
The Poincar\'e map $P$ from $\sum_0$ to itself is defined as:
\begin{equation*}
\begin{split}
P\colon\,\, &U\subset \sum\nolimits_0\longmapsto
\sum\nolimits_0\,\,,\\
&P=P_1^0\circ P_0^1\,\,.
\end{split}
\end{equation*}
\end{description}

\subsection{\subsectionfont Fixed Points of the Poincar\'e Map $P$}

On the Poincar\'e section $\sum_0$, we center the origin of
the coordinate frame at $w^{(\plus)}$, and denote the new
coordinates by $(x^{(0)}$, $z^{(0)}$, $v^{(\pm,0)})$.
On the Poincar\'e section $\sum_1$, we center the origin of
the coordinate frame at $w^{(-)}$, and denote the new
coordinates by $(x^{(1)},y^{(1)},z^{(1)},v^{(\pm,1)})$,
which satisfy the constraint equation
\begin{equation}
f^{(y)}\left(x^{(1)}, y^{(1)},z^{(1)}+\eta,
v^{(\plus,1)}+v_*^{\plus},v^{(-,1)}\right)=0\,\,
\label{constr}
\end{equation}
where for any $w\in\ScrH$,
$$
F^T(w)=\left(f^{(x)}(w),f^{(y)}(w),
f^{(z)}(w),f^{(v^+)}(w),f^{(v^-)}(w)\right)\,\,.
$$
Denote the vector $v^{\plus}$ in component form,
$v^{\plus}=(v_1^{\plus},v_2^{\plus},\dotsc,v_N^{\plus})'$.
Then we have the lemma,
\begin{lem}
$\frac{\partial f^{(y)}}{\partial z}(w^{(-)})$ and
$\frac{\partial f^{(y)}}{\partial v_j^+}(w^{(-)})$
($j=1,\dotsc,N$) cannot be zero simultaneously.
\end{lem}

\begin{proof}
Assume that they are zero simultaneously, then
$$
\tau\cdot\nabla f^{(y)}(w^{(-)})=0\,\,,
$$
where $\tau$ is the tangent vector of $h$ at $w^{(-)}$ and
``$\,\nabla\,$'' denotes gradient.
Then this implies that $h$ is tangent to $\sum_0$ at
$w^{(\plus)}$.
This contradiction proves the lemma.
\end{proof}

Let $\xi$ be one of the coordinates
$\{z,v_j^{\plus}\,(j=1,\dotsc,N)\}$, such that
$\frac{\partial f^{(y)}}{\partial\xi}(w^{(-)})\not=0$, and
denote by $v_{\plus}^{(1)}$ the vector
$(z^{(1)},v_j^{(\plus,1)}$ ($j=1,\dotsc,N))'\setminus
\{\xi\}$, i.e. with components consisting of
$\{z^{(1)},v_j^{(\plus,1)}$ ($j=1,\dotsc,N)\}$ without
$\xi$.

\begin{lem}
In a neighborhood of $w^{(-)}$, the Poincar\'e section
$\sum_1$ can be represented as a $C^2$ function
$$
\xi=\xi\left(x^{(1)},y^{(1)},v^{(-,1)},v_{\plus}^{(1)}\right)\,\,.
$$
\end{lem}

\begin{proof}
Applying the implicit function theorem to (\ref{constr}).
\end{proof}

\noindent
The map $P_0^1$ has the representation
$$
\begin{cases}
x^{(1)}=e^{-\alpha t_*}(x^{(0)}+x_*)\cos\,\beta\,t_*, &\\
y^{(1)}=e^{-\alpha t_*}(x^{(0)}+x_*)\sin\,\beta\,t_*, &\\
z^{(1)}+\eta=z^{(0)}e^{\gamma\,t_*}, &\\
v^{(\plus,1)}+v_*^{\plus}=e^{L^+ t_*}v^{(\plus,0)}, &\\
v^{(-,1)}=e^{L^{^{-}} t_*}(v^{(-,0)}+v_*^-)\,\,. &
\end{cases}
$$
The map $P_1^0$ can be approximated by its linearization at
$w^{(-)}$,
$$
\begin{pmatrix}
x^{(0)}\\
z^{(0)}\\
v^{(\plus,0)}\\
v^{(-,0)}
\end{pmatrix}
=A \begin{pmatrix}
x^{(1)}\\
y^{(1)}\\
z^{(1)}\\
v^{(\plus,1)}\\
v^{(-,1)}
\end{pmatrix} \ ,
$$
where
$$
A= \left ( \begin{array}{lcccr}
\frac{\partial f^{(x)}}{\partial x} & \frac{\partial f^{(x)}}
   {\partial y} & \frac{\partial f^{(x)}}{\partial z} &
   \frac{\partial f^{(x)}}{\partial v^{\plus}} &
  \frac{\partial f^{(x)}}{\partial v^-}\\ \\ 
\frac{\partial f^{(z)}}{\partial x} & \frac{\partial f^{(z)}}
  {\partial y} & \frac{\partial f^{(z)}}{\partial z} &
  \frac{\partial f^{(z)}}{\partial v^{\plus}} &
  \frac{\partial f^{(z)}}{\partial v^-}\\ \\ 
\frac{\partial f^{(v^{\plus})}}{\partial x} &
   \frac{\partial f^{(v^{\plus})}}{\partial y} &
   \frac{\partial f^{(v^{\plus})}}{\partial z} &
   \frac{\partial f^{(v^{\plus})}}{\partial v^{\plus}} &
   \frac{\partial f^{(v^{\plus})}}{\partial v^-}\\ \\ 
\frac{\partial f^{(v^-)}}{\partial x} &
   \frac{\partial f^{(v^-)}}{\partial y} &
   \frac{\partial f^{(v^-)}}{\partial z} &
   \frac{\partial f^{(v^-)}}{\partial v^{\plus}} &
   \frac{\partial f^{(v^-)}}{\partial v^-} \\
\end{array} \right )_{(w^{(-)})} \ .
$$
The constraint equation (\ref{constr}) can be approximated by
its linearization at $w^{(-)}$,
$$
B
\begin{pmatrix}
x^{(1)}\\
y^{(1)}\\
z^{(1)}\\
v^{(\plus,1)}\\
v^{(-,1)}
\end{pmatrix} = 0 \ ,
$$
where
$$ 
B = \left(\frac{\partial f^{(y)}}{\partial x} \
\frac{\partial f^{(y)}}{\partial y}\ \frac{\partial
f^{(y)}}{\partial z} \ \frac{\partial f^{(y)}}{\partial
v^{\plus}} \ 
\frac{\partial f^{(y)}}{\partial v^-}\right)_{(w^{(-)})}\ .
$$
With the above preparations, we can write the equations of
the fixed points of $P$ as follows in terms of the Silnikov
coordinates
$\left\{t_*,x^{(0)},v^{(-,0)},z^{(1)},v^{(\plus,1)}\right\}$
\begin{equation}
\begin{pmatrix}
x^{(0)}\hfill\\
\vspace{2pt}
(z^{(1)}+\eta)e^{-\gamma t_*}\hfill\\
\vspace{2pt}
e^{-L^+t_*}(v^{(\plus,1)}+
  v_*^{\plus})\hfill\\
\vspace{2pt}
v^{(-,0)}\hfill
\end{pmatrix}
=A\begin{pmatrix}
e^{-\alpha t_*}(x^{(0)}+x_*)\cos\,\beta\,t_*\\
\vspace{-3pt}
e^{-\alpha t_*}(x^{(0)}+x_*)\sin\,\beta\,t_*\\
\vspace{-3pt}
z^{(1)}\\
\vspace{-3pt}
v^{(\plus,1)}\\
\vspace{-3pt}
e^{L^{^{-}}t_*}(v^{(-,0)}+v_*^-)
\end{pmatrix}+\ScrR\,\,,
\label{fixe}
\end{equation}
where $\ScrR = \ScrO(e^{-\nu t_*})$ as $t_*\to+\infty$ for
some $\nu>0$.
By rescaling the coordinates as follows,
$$
 \left\{t_*,\xhat^{(0)}=x^{(0)}e^{\alpha t_*},
\vhat^{(-,0)}=v^{(-,0)}e^{\alpha t_*},
\zhat^{(1)}=z^{(1)}e^{\alpha t_*},\vhat^{(\plus,1)}=v^{(\plus,1)}
e^{\alpha t_*}\right\}\,\,,
$$
we can rewrite equation (\ref{fixe}) in the form 
\begin{equation}
\begin{pmatrix}
\xhat{^{(0)}}\\
\vspace{2pt}
0\\
\vspace{2pt}
0\\
\vspace{2pt}
\vhat{^{(-,0)}}
\end{pmatrix}
=A
\begin{pmatrix}
x_*\cos\beta t_*\\
\vspace{-3pt}
x_*\sin\beta t_*\\
\vspace{-3pt}
\zhat{^{(1)}}\\
\vspace{-3pt}
\vhat{^{(\plus,1)}}\\
\vspace{-3pt}
0
\end{pmatrix}+\ScrR_1\,\,,
\label{fixe1}
\end{equation}
where $\ScrR_1 = \ScrO(e^{-\nu_1t_*)}$ as $t_*\to +\infty$
for some $\nu_1>0$, and the constraint equation
(\ref{constr}) takes the form
\begin{equation}
B\begin{pmatrix}
x_*\cos\beta t_*\\
x_*\sin\beta t_*\\
\zhat{^{(1)}}\\
\vhat{^{(\plus,1)}}\\
0
\end{pmatrix}+\ScrR_2=0\,\,,
\label{constr1}
\end{equation}
where $\ScrR_2 = \ScrO(e^{-\nu_2t_*)}$ as $t_*\to +\infty$
for some $\nu_2>0$.
Solving the leading order term of (\ref{fixe1}) for
$(\zhat{^{(1)}},\vhat{^{(\plus,1)}})$, we have
\begin{equation}
C\begin{pmatrix}
\zhat{^{(1)}}\\
\vhat{^{(\plus,1)}}
\end{pmatrix}=-D
\begin{pmatrix}
x_*\cos\beta t_*\\
x_*\sin\beta t_*
\end{pmatrix}\,\,,
\label{lea1}'
\end{equation}
where
\[
C=\begin{pmatrix}
\frac{\partial f^{(z)}}{\partial z} &\frac{\partial f^{(z)}}
     {\partial v^{\plus}}\\
\frac{\partial f^{(v^+)}}{\partial z}
     &\frac{\partial f^{(v^+)}}{\partial v^{\plus}}
\end{pmatrix}_{(w^{(-)})},\,\, \quad 
D=\begin{pmatrix}
\frac{\partial f^{(z)}}{\partial x}
     &\frac{\partial f^{(z)}}{\partial y}\\
\frac{\partial f^{(v^+)}}{\partial x}
     &\frac{\partial f^{(v^+)}}{\partial y}
\end{pmatrix}_{(w^{(-)})}\,\,.\nonumber
\]

\begin{lem}
The matrix $C$ is invertible.
\end{lem}

\begin{proof}
Assume that $C$ is non-invertible; then there exists a
nonzero vector $\left(\begin{smallmatrix} \zhat{^{(1)}}\hfill\\
\vhat{^{(+,1)}}\hfill\end{smallmatrix}\right)$, such that
\[
C\begin{pmatrix}
\zhat{^{(1)}}\\
\vhat{^{(\plus,1)}}
\end{pmatrix}=0\,\,;
\]
thus
\[
w= A\begin{pmatrix}
0\\
0\\
\zhat{^{(1)}}\\
\vhat{^{(\plus,1)}}\\
0
\end{pmatrix}\in\ScrT_{w^{(\plus)}}W^s\cap\ScrT_{w^{(\plus)}}
 W^u
\]
where $\ScrT_{w^{(\plus)}}W^s$ is the tangent space of $W^s$
at $w^{(\plus)}$, similarly for $\ScrT_{w^{(\plus)}}W^u$.
Since $A$ is invertible, $w\not=0$.
We also know that
$$
\ScrT_{w^{(+)}}h\in\ScrT_{w^{(+)}}W^s\cap
\ScrT_{w^{(+)}}W^u\,\,.
$$
Since $\ScrT_{w^{(+)}}h$ is transversal to $\sum_0$ and
$w$ lies in $\sum_0$, we have
$$
\dim\left\{\ScrT_{w^{(+)}}W^u\cap\ScrT_{w^{(+)}}W^s
\right\}=2\,\,,
$$
which contradicts with Assumption (A4).
This completes the proof.
\end{proof}

\noindent
Solving (\ref{lea1}), we have
\begin{equation}
\begin{pmatrix}
\zhat{^{(1)}}\\
\vhat{^{(\plus,1)}}
\end{pmatrix}=-C^{-1}D
\begin{pmatrix}
x_*\cos\beta t_*\\
x_*\sin\beta t_*
\end{pmatrix}\,\,.
\label{lea2}
\end{equation}
Then solving the constraint equation (\ref{constr1}), to the
leading order, we have
$$
B\begin{pmatrix}
x_*\cos\beta t_*\hfill\\
x_*\sin\beta t_*\hfill\\
-C^{-1}D\begin{pmatrix}
x_*\cos\beta t_*\\
x_*\sin\beta t_*
\end{pmatrix}\\
\qquad 0\hfill
\end{pmatrix}=0\,\,,
$$
which can be rewritten as
\begin{equation}
\Delta_1\cos\beta t_*+\Delta_2\sin\beta t_*=0\,\,.
\label{lea3}
\end{equation}
If we assume condition 
\begin{description}
\item[\itemfont{(A5)}]
$\Delta_1$ and $\Delta_2$ do not vanish simultaneously,
\end{description}
then (\ref{lea3}) has a sequence of solutions:
\begin{equation}
t_*^{(\ell)}=\frac1\beta [\ell\pi-\varphi],\,\,
  \ell\in Z^{\plus},\label{lea4}
\end{equation}
where
\[
\varphi=\arctan\left\{\Delta_1/\Delta_2\right\}\,\,.
\]
Then substituting (\ref{lea2}) and (\ref{lea4}) into the
leading order terms of (\ref{fixe1}), we can solve for
$\left(\begin{smallmatrix}
\xhat{^{(0)}}\hfill\\
\vhat{^{(-,0)}}\hfill\end{smallmatrix}\right)$ to the
leading order.
Finally, applying the implicit function theorem, we have the
fixed point theorem.

\begin{thm}
[Fixed Point Theorem] Under the assumption (A5), there exists an integer
$\ell_0>0$, such that there exists a sequence of solutions
to the equations (\ref{fixe1}) and (\ref{constr1}) labeled
by $\ell$ ($\ell\Ge \ell_0$):
\begin{eqnarray*}
& & t_*=T^{(\ell)}\,, \quad \xhat{^{(0)}}=x^{(\ell)},\,\,
  \vhat{^{(-,0)}}=v_{(-,0)}^{(\ell)}\\
& & \zhat{^{(1)}}=z^{(\ell)}\,, \quad \vhat{^{(+,1)}}=
  v_{(\plus,1)}^{(\ell)}\,\,;
\end{eqnarray*}
where, as $t\to+\infty$,
$$
T^{(\ell)}=\frac1\beta\,\,[\ell\pi-\varphi]+o(1)\,\,.
$$
\end{thm}
\noindent
For a complete proof of this theorem, see \cite{Li99}.

\subsection{\subsectionfont Smale Horseshoes}
\begin{description}
\item[\itemfont{Definition.}]
For sufficiently large number $\ell$, we define the slab
$S_\ell$ in $\sum\nolimits_0$ as follows:
\begin{eqnarray*}
S_\ell &=& \Biggl\{w\in\sum\nolimits_0 \Bigm\arrowvert
\eta\exp \left\{-\gamma\left(T^{(2(\ell+1))}-\frac{\pi}{2\beta}
     \right)\right\}\Le z^{(0)}(w) \\
& & \Le \eta\exp\left\{-\gamma\left(T^{(2\ell)}-\frac{\pi}{2\beta}
     \right)\right\}\,\,,\\
& & \vert x^{(0)}(w)\vert\Le\eta\exp\left\{-\frac12\,
     \alpha T^{(2\ell)}\right\}\,\,,\\
& & \| v^{(-,1)}(P_0^1(w))\| \Le\eta\exp\left\{-\frac12\,
     \alpha T^{(2\ell)}\right\}\,\,,\\
& & \| v^{(\plus,1)}(P_0^1(w))\| \Le \eta\exp\left\{
     -\frac12\,\alpha T^{(2\ell)}\right\}\Biggr\}\,\,,
\end{eqnarray*}
so that it contains two fixed points of $P$ denoted by
$p_\ell^{\plus}$ and $p_\ell^-$.
\end{description}

\begin{figure}
{\begin{picture}(6,2)
\put(.75,0){\vector(0,1){2}}
\put(0,1){\vector(1,0){1.5}}
\put(.25,1.25){\line(0,1){.25}}
\put(.25,1.25){\line(1,0){1}}
\put(.25,1.5){\line(1,0){1}}
\put(1.25,1.25){\line(0,1){.25}}
\put(.1,1.35){$2$}
\put(1.35,1.35){$1$}
\put(.6,1.1){$4$}
\put(.6,1.6){$3$}
\put(.85,1.875){$e_{z^{(0)}}$}
\put(1.6,1){$e_{x^{(0)}}$}
\put(2.75,.75){\includegraphics{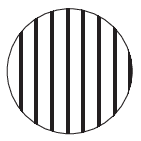}}
\put(2.8,1.4){$\mathbf{e}_{v^{(+,0)}}$}
\put(4,.25){\includegraphics{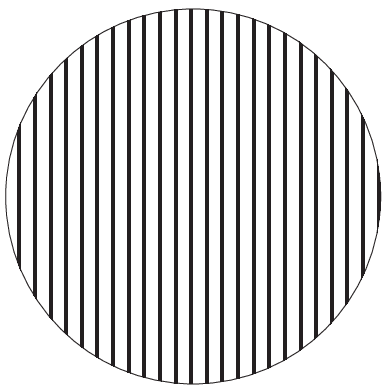}}
\put(4.5,1.85){$\mathbf{e}_{v^{(-,0)}}$}
\put(2.25,1.5){$\otimes$}
\put(3.65,1.5){$\otimes$}
\end{picture}}
\caption{The product representation of $S_\ell$\ .}
\label{prei}
\end{figure}

\begin{figure}{\begin{picture}(6,2)
\put(0,0){\includegraphics{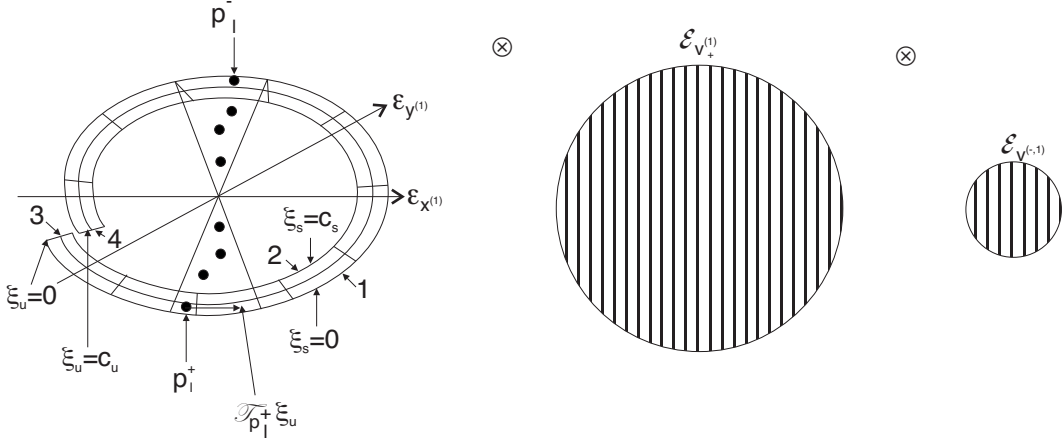}}
\end{picture}}
\caption{A geometric representation of $\nabla P_1^0(w^{(-)})\circ
\nabla P_0^1(w^{(\plus)})\circ S_\ell$\ .}
\label{image}
\end{figure}

We choose a basis on the tangent space
$\ScrT_{w^{(-)}}\sum_1$ represented in the coordinates 
$$(x^{(1)}, y^{(1)}, v^{(-,1)}, \xi, v^{(1)}_+)$$ as follows:
\begin{eqnarray*}
& & E_{x^{(1)}} =\left(1,0,0,\nabla\xi(w^{(-)})\circ
     (1,0,0,0),0\right)\,\,,\\
& & E_{y^{(1)}} =\left(0,1,0,\nabla\xi(w^{(-)})\circ
     (0,1,0,0),0\right)\,\,,\\
& & \bfE_{v^{(-,1)}} =\left(0,0,\one,\nabla\xi(w^{(-)})\circ
     (0,0,\one,0),0\right)\,\,,\\
& & \bfE_{v_{\plus}^{(1)}} =\left(0,0,0,\nabla\xi(w^{(-)})\circ
     (0,0,0,\one),\one\right)\,\,;
\end{eqnarray*}
where $\one$ represents a basis for the corresponding
components.
Denote by $\{e_{_{x^{(0)}}},e_{_{z^{(0)}}},{\bf e}_{_{v^{(+,0)}}}$,\break
${\bf e}_{_{v^{(-,0)}}}\}$ the unit vectors along
$(x^{(0)},z^{(0)},v^{(+,0)},v^{(-,0)})$-directions in
$\sum_0$.
In this coordinate frame, $S_\ell$ has the product
representation as shown in Fig.\ref{prei}.
\noindent
Under the linear map $A$ (the linearization of $P_1^0$ at
$w^{(-)}$), the coordinate frame
$\left\{E_{x^{(1)}},
E_{y^{(1)}},\bfE_{v^{(-,1)}},\bfE_{v_+^{(1)}}\right\}$
is mapped into a coordinate frame
$\left\{\ScrE_{x^{(1)}},\ScrE_{y^{(1)}},
\bfscrE_{v^{(-,1)}},\bfscrE_{v_+^{(1)}}\right\}$ on $\sum_0$ with
origin at $w^{(\plus)}$.
In this coordinate frame, $\nabla P_1^0(w^{(-)})\circ
\nabla P_0^1(w^{(\plus)})\circ S_\ell$ has the
representation as shown in Fig.\ref{image} on $\sum_0$.
\noindent
We introduce a system of curvilinear coordinates
$(\xi_u,\xi_s)$ on the
$(\ScrE_{x^{(1)}},\ScrE_{y^{(1)}})$-plane such that
$$
\{\xi_u=0\},\,\,
\{\xi_u=c_u\text{ (a constant)}\},\,\,
\{\xi_s=0\},\,\,
\{\xi_s=c_s\text{ (a constant)}\}
$$
correspond to the boundaries $3$, $4$, $1$, $2$ of the
annulus on the $(\ScrE_{x^{(1)}},\ScrE_{y^{(1)}})$-plane.
Let $\ScrE_\xi^{\plus}=\ScrT_{p_\ell^+}\xi_u$ be the tangent
vector to the $\xi_u$ coordinate at $p_\ell^{\plus}$.
We make the following assumption:
\begin{description}
\item[\itemfont{(A6)}]
$\ \ \ \ \textrm{Span}\,\{e_{_{x^{(0)}}},\bfe_{_{v^{(-,0)}}},
\ScrE_\xi^{\plus},\bfscrE_{v_+^{(1)}}\}=\sum_0\ $.
\end{description}

Under the assumptions (A1) -- (A6), we can verify the
Conley-Moser conditions \cite{Mos73} in the same way as in \cite{Li99},
which lead to the construction of Smale horseshoes.
Let $\ScrW$ be a set which consists of elements of the
doubly infinite sequence form 
$$
a=(\,\cdots a_{-2}a_{-1}a_0,\,a_1a_2\cdots\,)\,\,,
$$
where $a_k\in\{0,1\}$, $k\in Z$.
We introduce a topology in $\ScrW$ by taking as neighborhood
basis of
\[
a^*=(\,\cdots a_{-2}^*a_{-1}^*a_0^*,\,a_1^*a_2^*\cdots\,)\,\,,
\]
the set
\[
W_j=\left\{a\in\ScrW\bigm\arrowvert a_k=a_k^*\,\,(\vert
k\vert<j)\right\}
\]
for $j=1,2,\ldots\,\,$.
This makes $\ScrW$ a topological space.
The shift automorphism $\chi$ is defined on $\ScrW$ by
\begin{eqnarray*}
\chi &\colon & \ScrW\longmapsto\ScrW\,\,,\\
& & \forall\,a\in\ScrW,\,\,
\chi(a)=b,\quad\text{where}\quad
b_k=a_{k+1}\,\,.
\end{eqnarray*}
The shift automorphism $\chi$ exhibits sensitive dependence 
on initial conditions, which is the {\em{hallmark}} of {\em{chaos}}.
\begin{thm}
[Smale Horseshoe Theorem]Under the assumptions (A1) -- (A6) 
for the evolution equation (\ref{ee}), for all sufficiently large integers 
$\ell$, there exists a sequence of compact Cantor subsets $\Lambda_\ell$ 
of $S_\ell$, $\Lambda_\ell$ consists of points and is invariant under
$P$. $P$ restricted to $\Lambda_\ell$ is topologically conjugate
to the shift automorphism $\chi$ on two symbols $0$ and
$1$. That is, there exists a homeomorphism
$$
\phi_\ell\colon\, \ScrW\longmapsto \Lambda_\ell\,\,,
$$
such that the following diagram commutes:
$$
\begin{CD}
\ScrW @>{\phi_\ell}>> \Lambda_\ell\\
@V{\chi}VV @VV{P}V\\
\ScrW @>>{\phi_\ell}> \Lambda_\ell\,\,.
\end{CD}
$$
\end{thm}

\begin{proof}
With preliminaries given above, the proof follows in the same
way as in \cite{Li99}.
\end{proof}

\section{Conclusion and Discussion}

In this note, we have generalized the construction of Smale
horseshoes in \cite{Li99} to a general evolution equation,
which indicates that our techniques \cite{LW97}, \cite{Li99}
on constructing Smale horseshoes in a neighborhood of a
Silnikov homoclinic orbit has a much wider application.
On the other hand, so far we can only handle evolution
equations with finitely many linearly unstable modes.
For evolution equations with infinitely many linearly
unstable modes, we cannot invert certain linear operators in
establishing the existence of fixed points of the Poincar\'e
map.
Nevertheless, this note furnishes an initiation for studying
Silnikov homoclinic orbits for general evolution equations,
thereby proving the existence of chaos for general evolution
equations.

\end{document}